\documentclass{amsart}
\usepackage{amscd,amssymb,hyperref}

\theoremstyle{plain}
\newtheorem{prop}[subsection]{Proposition}
\newtheorem{thm}[subsection]{Theorem}
\newtheorem*{theo}{Theorem}

\newtheorem{lem}[subsection]{Lemma}
\newtheorem{cor}[subsection]{Corollary}
\newtheorem{definition}[subsection]{Definition}
\theoremstyle{remark}
\newtheorem{rem}[subsection]{Remark}

\theoremstyle{definition}

\numberwithin{equation}{section}

\renewcommand{\b}[1]{\mathbf{#1}}

\newcommand{\A}{{\mathcal A}}
\newcommand{\Ai}{{\mathcal A}_\infty}

\newcommand{\cC}{{\mathcal C}}

\newcommand{\cE}{{\mathcal E}}
\newcommand{\cF}{{\mathcal F}}
\newcommand{\cG}{{\mathcal G}}
\newcommand{\LL}{{\mathcal L}}
\newcommand{\cR}{{\mathcal R}}
\newcommand{\cS}{{\mathcal S}}
\newcommand{\cT}{{\mathcal T}}

\newcommand{\Z}{{\mathbb Z}}

\newcommand{\cn}{{\mathbb C}}
\newcommand{\R}{{\mathbb R}}
\newcommand{\C}{{\mathbb C}}

\newcommand{\CP}{{\mathbb{CP}}}
\newcommand{\cp}{\cn{\mathbb P}}
\newcommand{\bS}{{\mathbb S}}

\newcommand{\bl}{{\boldsymbol{\lambda}}}
\newcommand{\by}{{\boldsymbol{y}}}
\newcommand{\sA}{{\mathsf A}}

\newcommand{\sfB}{{\sf B}}
\newcommand{\sfD}{{\mathsf D}}
\newcommand{\sfb}{{\sf b}}

\newcommand{\sfI}{{\mathsf I}}

\newcommand{\sK}{{\mathsf K}}
\newcommand{\sM}{{\sf M}}
\newcommand{\sfM}{{\sf M}}
\newcommand{\sfX}{{\mathsf X}}
\newcommand{\sfY}{{\mathsf Y}}

\newcommand{\al}{{\alpha }}
\newcommand{\la}{{\lambda }}
\newcommand{\nbc}{{\mathbf{nbc}}}

\newcommand{\bul}{{\bullet }}
\newcommand{\p}{\partial}
\newcommand{\ra}{\rightarrow}
\newcommand{\bone}{{\mathbf 1}}
\renewcommand{\c}{{\gamma }}
\renewcommand{\ll}{{\ell }}
\newcommand{\D}{{\Delta}}
\newcommand{\TT}{{({\mathbb C}^*)^n}}

\renewcommand{\L}{{\Lambda }}

\DeclareMathOperator{\rank}{rank}
\DeclareMathOperator{\codim}{codim}

\DeclareMathOperator{\ii}{i}

\DeclareMathOperator{\End}{End}
\DeclareMathOperator{\Aut}{Aut}
\DeclareMathOperator{\Dep}{Dep}

\DeclareMathOperator{\Mat}{Mat}
\DeclareMathOperator{\Hom}{Hom} 
\DeclareMathOperator{\GL}{GL} 

\begin{document}

\title[Stratified Morse Theory in Arrangements]
{Stratified Morse Theory in Arrangements}
\author[D.~Cohen]{Daniel C.~Cohen$^\dag$}
\address{Department of Mathematics, Louisiana State University,
Baton Rouge, LA 70803}
\email{\href{mailto:cohen@math.lsu.edu}{cohen@math.lsu.edu}}
\urladdr{\href{http://www.math.lsu.edu/~cohen/}
{http://www.math.lsu.edu/\~{}cohen}}
\thanks{{$^\dag$}Partially supported by National Security Agency grant
H98230-05-1-0055}

\author[P.~Orlik]{Peter Orlik}
\address{Department of Mathematics, University of Wisconsin,
Madison, WI 53706}
\email{\href{mailto:orlik@math.wisc.edu}{orlik@math.wisc.edu}}

\dedicatory{For Robert MacPherson on the occasion of his sixtieth birthday.}



\begin{abstract}
This paper is a survey of our work based on the stratified Morse theory
of Goresky and MacPherson.  First we discuss the Morse theory 
of Euclidean space stratified by an arrangement.  This is used to show that 
the complement of a complex hyperplane arrangement admits a minimal 
cell decomposition.  Next we review the construction of 
a cochain complex whose cohomology computes the 
local system cohomology
of the complement of a complex hyperplane arrangement.
Then we present results on the Gauss-Manin connection for the
 moduli space of arrangements of a fixed combinatorial type in
 rank one local system cohomology.
\end{abstract}


\maketitle

\section{Introduction}
\label{sec:intro}

Let $V$ be a  complex vector space of dimension $\ell\geq 1$. 
A hyperplane arrangement $\A=\{H_1,\ldots,H_n\}$  is a set of $n \geq 0$
hyperplanes in $V$. Let $\sfM=V\setminus \bigcup_{j=1}^n H_j$ denote the complement. 
Introduce coordinates $u_1,\ldots,
u_{\ell}$ in $V$ and, for each $j$, $1\le j\le n$, 
choose a degree one polynomial $\al_j$  so that the 
hyperplane $H_j \in \A$ is defined by the vanishing of $\al_j$. 
Let $\bl=(\la_1, \ldots,\la_n)$ be a set of complex weights for
the hyperplanes. 
Given $\bl$,
we define a multivalued holomorphic function on $\sfM$ by
\begin{equation*}
\Phi(u;\bl)=\prod_{j=1}^n\al_j^{\la_j}.
\end{equation*}

A generalized hypergeometric integral is of the form
\[
\int_{\sigma}\Phi(u;\bl) \eta
\] 
where $\sigma$ is a suitable domain
of integration and $\eta$ is a holomorphic form on $\sfM$, see 
\cite{AoK1}. When $\ell=1$, $n=3$ and $\al_1=u, \al_2=u-1, \al_3=u-x$, 
this is the Gauss hypergeometric integral.   
Selberg's
integral \cite{Sel1} is another special case:
\[
\int_0^1 \cdots \int_0^1 (u_1  \cdots u_{\ell})^{x-1}
[(1-u_1)\cdots(1-u_{\ell})]^{y-1} |\Delta(u)|^{2z}\, du_1 \ldots
du_{\ell}
\]
where $\Delta(u)=\prod_{i<j}(u_j-u_i)$.
Hypergeometric integrals occur in the representation theory of Lie
algebras and quantum groups \cite{SV,Var2}. In physics, these hypergeometric
integrals form solutions to the Knizhnik-Zamolodchikov
differential equations in conformal field theory \cite{ScV1,Var2}. The
space of integrals is identified with
a cohomology group, $H^\ell(\sfM;\LL)$, of the complement, 
with coefficients in a complex rank one local system.  
Associated to $\bl$, there is
a rank one representation $\rho:\pi_1(\sfM) \to \cn^*$, given by
$\rho(\c_j) =t_j$, where $\b{t}=(t_1,\dots,t_n)\in (\cn^*)^n$ is
defined by $t_j= \exp(-2\pi\ii\la_j)$, and $\c_j$ is any meridian
loop about the hyperplane $H_j$ of $\A$, and a corresponding rank
one local system $\LL=\LL_{\b{t}}=\LL_{\bl}$ on $\sM$.  Equivalently, 
weights $\bl$ determine a flat connection on the trivial line bundle over 
$\sfM$, with connection one-form $\omega_\bl = d\log \Phi(u;\bl)$.

The first problem is to calculate the local system cohomology groups
$H^q(\sM,\LL)$. 
The methods used by Aomoto and Kita \cite{AoK1}, Esnault, Schechtman, and Viehweg
\cite{ESV}, Schechtman, Terao, and Varchenko \cite{STV} and others are 
described in detail in \cite{OrT2}.
These use the twisted de Rham
complex, $(\Omega^\bul(*\A),\nabla)$, of global rational differential forms on $V$ with
arbitrary poles along the divisor $\bigcup_{j=1}^n H_j$, with differential $\nabla(\eta) = d\eta+\omega_\bl \wedge \eta$.  The cochain
groups of this complex are infinite dimensional. Conditions must be imposed on
the weights in order to reduce the problem to a finite dimensional
setting. These are the nonresonance conditions of \cite{STV}.  
Under these conditions, the calculation may be 
reduced to combinatorics and yields 
$H^q(\sfM;\LL) = 0$ for $q\neq\ll$ and
$\dim H^\ll(\sfM;\LL) = |e(\sfM)|$, where $e(\sfM)$ is the Euler characteristic
of the complement. This approach provides less 
information for resonant weights, those for which the aformentioned 
nonresonance conditions do not hold. By contrast, the 
results obtained below using stratified Morse theory are valid for
arbitrary weights. 

Weights $\bl$ give rise to a local system on the complement of every arrangement 
that is combinatorially equivalent to $\A$.  The resulting local system cohomology groups 
comprise a flat vector bundle over the moduli space of such arrangements.  
The second problem is to determine the Gauss-Manin
connection in this cohomology bundle.  
For instance, the Gauss hypergeometric function is defined on the complement of
the arrangement of three points in $\cn$. It satisfies a second
order differential equation which, when converted into a system of
two linear differential equations, may be interpreted as a
Gauss-Manin connection on the moduli space of arrangements of the
same combinatorial type \cite{OrT2}. This idea has been
generalized to all arrangements by Aomoto \cite{Aom1} and Gelfand
\cite{Gel1}. The connection is obtained by differentiating in the
moduli space. For arrangements in general position, and nonresonant weights, 
explicit connection matrices were
obtained by Aomoto and Kita \cite{AoK1}.  Unfortunately,
this pioneering work is available only in Japanese. The relevant
matrices have been reproduced in \cite{OrT2,CO4}. 

The (flat) Gauss-Manin connection in the cohomology bundle corresponds 
to a representation of the fundamental group of the moduli space.  The endomorphisms 
arising in the connection one-form, which we refer to as Gauss-Manin endomorphisms, 
may be realized as logarithms of certain automorphisms.  This interpretation, used in \cite{CO5}, 
allows for local calculations, valid for all arrangements and all weights.  This paper is a survey 
of our work on these problems. 

Section \ref{sec:smt} presents basic results on the Morse theory of 
Euclidean space stratified by an arrangement (of subspaces), following \cite{GM,C1}.  
In Section \ref{sec:min}, we use these results to give a proof of 
a theorem of Dimca and Papadima \cite{DP1} and Randell \cite{Ra2}, 
which asserts that the complement of a complex hyperplane arrangement 
is a minimal space, admiting a cell decomposition for which the 
number of $q$-cells is equal to the $q$-th Betti number for each $q$.  
In Section \ref{sec:local systems}, we review
the stratified Morse theory construction from \cite{C1,CO1} of a finite
cochain complex $(K^\bul(\A),\D^\bul)$, the cohomology of which is
naturally isomorphic to $H^*(\sfM;\LL)$.  This leads to the
construction of the universal complex $(\sK^\bul,\D^\bul(\b{x}))$
for local system cohomology, where $\sK^q=\Lambda\otimes_\cn K^q$ and
$\Lambda=\cn[x_1^{\pm 1},\dots,x_n^{\pm 1}]$. 
We recall a combinatorial model
for $H^*(\sfM,\cn)$, called
the Orlik-Solomon algebra, $A(\A)$.  
The one-form $\omega_\bl$ corresponds to an element $a_\bl$ of the 
Orlik-Solomon algebra.  Multiplication by $a_\bl$ gives this algebra 
the structure of a cochain complex, $(A^\bul(\A),a_\bl)$.  
The Aomoto complex is the universal complex for this 
cochain complex.     It is
chain equivalent to the linearization of the universal complex. 
This informs on the relationship between the characteristic varieties 
of complements of arrangements (jumping loci for local system cohomology) 
and the resonance varieties of arrangements (jumping loci for 
the cohomology of the Orlik-Solomon complex).

In Section \ref{sec:moduli}, we move from consideration of a fixed 
arrangement to the study of all arrangements of a given combinatorial
type. We define  the moduli space of
arrangements with a fixed combinatorial type $\cT$ and the set $\Dep(\cT)$ of
dependent sets in type $\cT$.
We present results concerning the homology of the moduli space.

In Section \ref{sec:GM}, we work with
 a smooth, connected component of the moduli space, $\sfB(\cT)$.
There is a fiber bundle $\pi:\sfM(\cT) \to \sfB(\cT)$ whose
fibers, $\pi^{-1}(\sfb)=\sfM_{\sfb}$, are complements of
arrangements $\A_{\sfb}$ of type $\cT$.  Since $\sfB(\cT)$ is
connected, $\sfM_{\sfb}$ is diffeomorphic to $\sfM$. The fiber
bundle $\pi:\sfM(\cT) \to \sfB(\cT)$ is locally trivial.
Given a local system on the fiber, consider the
associated flat vector bundle
$\pi^{q}:\b{H}^{q}(\LL)\to\sfB(\cT)$, with fiber
$(\pi^{q})^{-1}(\sfb)=H^{q}(\sfM_{\sfb};\LL_{\sfb})$ at
$\sfb\in\sfB(\cT)$ for each $q$, $0\le q \le \ell$. Fixing a
basepoint $\sfb \in\sfB(\cT)$, the operation of parallel
translation of fibers over curves in $\sfB(\cT)$ in the vector
bundle $\pi^{q}:\b{H}^{q}(\LL)\to\sfB(\cT)$ provides a complex
representation 
\[
\Psi^{q}_{\cT}:\pi_{1}(\sfB(\cT),\sfb)
\longrightarrow \Aut_\cn(H^{q}(\sfM_{\sfb};\LL_{\sfb})).
\]  
The
loops of primary interest are those linking moduli spaces of
codimension one degenerations of $\cT$.  Such a degeneration is a
type $\cT'$ whose moduli space $\sfB(\cT')$ has codimension one in
the closure of $\sfB(\cT)$. In this case we say 
that $\cT$ covers $\cT'$.

When $\cT$ covers $\cT'$ and
$\c\in\pi_1(\sfB(\cT),\sfb)$ is a simple loop linking $\sfB(\cT')$
in $\overline{\sfB(\cT)}$, write $\Psi^q_\cT(\c) = \exp(-2\pi\ii \Omega)$.  
We denote this Gauss-Manin
endomorphism in the bundle $\pi^q:\b{H}^q(\LL) \to \sfB(\cT)$ by
$\Omega=\Omega^q_{\LL}(\sfB(\cT'),\sfB(\cT))$. The rest of this survey reports on 
our results concerning these endomorphisms.

In Section \ref{sec:formal}, for each subset $S$ of hyperplanes, we define an
endomorphism $\tilde{\omega}^\bul_S$ of the Aomoto complex of
a general position arrangement of $n$ hyperplanes in $\cn^\ell$.
When $\cT$ covers $\cT'$, we construct a suitable linear combination 
of these maps, which induces an endomorphism of the Aomoto complex of type $\cT$.  
The specialization $\b{y} \mapsto \bl$ in the Aomoto complex then yields an endomorphism
$\omega^\bul_{\bl}(\cT',\cT)$ of the 
Orlik-Solomon complex 
$A^\bul(\cT) = A^\bul(\A)$ of an arrangement $\A$ 
of type $\cT$. This leads to our main result, stated in more detail in Section \ref{sec:formal}.

\begin{theo}[\cite{CO5}]
Let $\sfM$ be the complement of an arrangement $\A$ of
type $\cT$ and let $\LL$ be the local system on $\sfM$ defined by
weights $\bl$. Suppose $\cT$ covers $\cT'$.
Then there is a surjection $\xi^q:A^q(\cT) \twoheadrightarrow H^q(\sfM,\LL)$
so that the Gauss-Manin endomorphism
$\Omega^q_{\LL}(\sfB(\cT'),\sfB(\cT))$ in local system cohomology
is determined by the equation
\[
\xi^q \circ \omega^q_{\bl}(\cT',\cT) =
\Omega^q_{\LL}(\sfB(\cT'),\sfB(\cT)) \circ \xi^q. \qed
\]
\end{theo}

In Section \ref{sec:eigen}, we report on the spectrum of
the Gauss-Manin endomorphism.  The pair
$(\cT',\cT)$ determines a set of hyperplanes $S\subset \A$ 
and an integer $r$. We call $(S,r)$  the principal dependence. Let
$\la_S=\sum_{H_j\in S}\la_j$. 

\begin{theo}[\cite{CO6}]
Suppose $\cT$ covers $\cT'$ with principal dependence $(S,r)$. 
Let $\bl$ be a collection of weights satisfying $\la_S \neq 0$.  
Then the Gauss-Manin endomorphism  $\Omega^q_{\LL}(\sfB(\cT'),\sfB(\cT))$
is diagonalizable, with spectrum contained in $\{0,\la_S\}$. \qed
\end{theo} 

We illustrate these results with an example in Section \ref{sec:selberg}.

\section{Morse Functions for Arrangements}
\label{sec:smt}

Goresky and MacPherson developed stratified Morse theory in order to extend the 
class of spaces to which Morse theory applies.  This generalization may be used to 
study singular spaces, noncompact spaces, etc.  The latter is illustrated in Part III of 
their book \cite{GM} using real subspace arrangements.  The topology of the complement 
of such an arrangement is analyzed by Morse theoretic means, by considering the 
stratification of the ambient Euclidean space determined by the arrangement and realizing 
the complement as one of the strata. We recall some of their constructions and results needed in 
the sequel.

Let $V$ be a real vector space of dimension $\ell \geq 1$, 
and let $\A$ 
be an arrangement of affine
subspaces in $V$.  An edge (or flat) of $\A$ is a nonempty 
intersection $X$ of elements of $\A$.  Let $L=L(\A)$ be the 
set of all edges of $\A$.  Unless otherwise noted, we partially order the set $L$ by reverse inclusion.

The arrangement $\A$ gives rise to a Whitney stratification 
$\cS$ of $V$ with a stratum
\[
\cS_X  = X \setminus \bigcup_{Y \subsetneq X} Y
\]
for each edge $X \in L$.  The complement 
$\sfM$ of $\A$ 
is the stratum corresponding to the edge
$V$ (the intersection of no elements of $\A$). For any 
edge $X$, the closure of $\cS_X$ is $X$.
Note that a complex hyperplane arrangement may be viewed 
as a real subspace arrangement with even-dimensional strata.

For almost any point $p \in \sfM$, the function $f:V \to \R$
given by 
\begin{equation} \label{eq:dist}
f(u) = [\operatorname{distance}(p,u)]^2
\end{equation}
is a Morse function on $V$ with respect to the stratification 
$\cS$, see \cite[I.2.2]{GM}.  For $r\in\R$, let
\[
\sfM_{\le r} = \{u \in \sfM \mid f(u) \le r\}.
\]
The function $f$ has a unique 
critical point  on each edge. It is  a minimum.  Furthermore, 
 Goresky and MacPherson 
show in \cite[III.3]{GM} that the Morse function $f$ is \emph{perfect}:  
if $v\in \R$ is a critical value and $\epsilon>0$ is 
sufficiently small, the long exact homology 
sequence of the pair $(\sfM_{\le v+\epsilon},\sfM_{\le v-\epsilon})$ 
splits into short exact sequences
\begin{equation} \label{eq:perfect}
0 \to H_q(\sfM_{\le v-\epsilon};\Z) \to H_q(\sfM_{\le v+\epsilon};\Z) 
\to H_q(\sfM_{\le v+\epsilon},\sfM_{\le v-\epsilon};\Z) \to 0.
\end{equation}
Using this, they calculate the homology $H_*(\sfM;\Z)$ in terms 
of the poset $L$ (ordered by inclusion), 
see \cite[III.1.3. Theorem A]{GM}.  This result has prompted 
a great deal of work on the cohomology of the complement of a 
subspace arrangement, culminating with the determination of 
the cup product structure of this cohomology ring by de Longeville and Schultz \cite{dLS} 
and Deligne, Goresky, and MacPherson \cite{DGM}.

One can produce a Morse function such as \eqref{eq:dist} that meets the strata of $V$ according to codimension.

\begin{definition} \label{def:wsi}
Let $Z$ be a Whitney stratified subset of Euclidean space.  A Morse function $f:Z \to \R$ is said to be \emph{weakly self-indexing} with respect to the stratification $\{S_\alpha\}$ of $Z$ if for each $q$, $0 \le q \le \dim Z$, we have
\[
\max_{\codim S_\alpha = q-1}\{\text{critical values of }f \mid S_\alpha\}
<
\min_{\codim S_\beta = q}\{\text{critical values of }f \mid S_\beta\}.
\]
\end{definition}

\begin{prop} \label{prop:wsi}
Let $\A$ be an arrangement of subspaces in the real vector space $V$.  Then there is a positive definite quadratic form $f:V \to \R$ which is a weakly self-indexing Morse function with respect to the stratification $\{S_\alpha\}$ of $V$ given by $\A$, whose critical points consist of a unique minimum on each stratum. \qed
\end{prop}

The proof of this result given in \cite[\S1]{C1} shows that there are choices of coordinates $\{u_i\}$ on $V$ and positive constants $\{\omega_i\}$, $1\le i \le \ell=\dim V$, for which the quadratic form 
$f(u_1,u_2,\dots,u_\ell) = \sum_{i=1}^\ell \omega_i u_i^2$
is a weakly self-indexing Morse function with respect to the stratification determined by $\A$.  This provides an inductive algorithm for the construction of a complete flag in $V$ that is transverse to the arrangement $\A$. 

In the rest of this paper, we return to the special case of a complex hyperplane arrangement 
where we can say more.

\section{Minimality}
\label{sec:min}

The notion of \emph{minimality} has played a significant role in recent work on the topology of arrangements, see for instance the work of
Papadima and Suciu \cite{PS}, Dimca and Papadima \cite{DP1,DP2}, and Randell \cite{Ra2}.

\begin{definition} \label{def:minimal}
A space $\sfX$ is said to be \emph{minimal} if $\sfX$ has the homotopy type of a connected, finite-type CW-complex $W$ such that, for each $q\ge 0$, the number of $q$-cells in $W$ is equal to the rank of $H_q(\sfX;\Z)$.
\end{definition}
Note that, for a minimal space $\sfX$, all homology groups $H_q(\sfX;\Z)$ are finitely generated and torsion-free.  If $\sfX$ is a $1$-connected space with the homotopy type of a connected, finite-type CW-complex, and the homology of $\sfX$ is torsion-free, then $\sfX$ is minimal by work of Anick \cite{An}. However, many spaces (with torsion-free homology) are not minimal.  For instance, the complement of a non-trivial knot does not admit a minimal cell decomposition.

Dimca and Papadima \cite{DP1} and Randell \cite{Ra2} used various forms of Morse theory to show that the complement of a complex hyperplane arrangement is minimal.  This result may also be established using stratified Morse theory.

\begin{thm} \label{thm:min}
 Let $\A=\{H_1,\dots,H_n\}$ be a complex hyperplane arrangement in the complex vector space $V \cong \C^\ell$. Then the complement $\sfM=\sfM(\A) = V \setminus \bigcup_{i=1}^n H_i$
 is a minimal space.
\end{thm}
\begin{proof}
Without loss of generality, assume that $\A$ is an \emph{essential} arrangement in $\C^\ell$, that is, 
$\A$ contains $\ell$ linearly independent hyperplanes.  Then the edges of $\A$ have codimensions $0$ through $\ell$.  The proof is by induction on $\ell$.

In the case $\ell=1$, $\A$ is a finite collection of points in $V=\C$, 
and the complement $\sfM(\A)$ has the homotopy type of a bouquet of circles, which is a minimal space.

For general $\ell$, let 
\begin{equation} \label{eq:flag}
\cF:\quad \emptyset = \cF^{-1} \subset \cF^0 \subset \cF^1 \subset
\cF^2 \subset\dots \subset \cF^{\ll-1} \subset \cF^\ll = V,
\end{equation}
be a complete flag 
in $V=\C^\ell$ that is transverse to the Whitney stratification of $V$ determined by $\A$.  Choose coordinates $\{u_i\}$ so that, for each $k\le\ell-1$,  $\cF^{k}=\{u_{k+1}=\dots=u_\ell=0\}$.  Let $f:V\to\R$ be a Morse function ``about'' the flag $\cF$ that is weakly self-indexing with respect to the stratification of $V$ determined by $\A$. 

Since $f$ is weakly self-indexing, there are constants $a$ and $b$ so that all critical values of $f$ on edges of codimension less than $\ell$ are smaller than $a$, and all critical values of $f$ on edges of codimension $\ell$ are in the interval $(a,b)$.  
For such $a$ and $b$, 
$\sfM_{\le b}$ is a deformation retract of the complement $\sfM$ of $\A$, and $\sfM \cap \cF^{\ell-1}$ is a deformation retract of $\sfM_{\le a}$.  Since $\sfM \cap \cF^{\ell-1}$ is the complement of the arrangement $\A\cap\cF^{\ell-1}$ in $\cF^{\ell-1}=\C^{\ell-1}$, by induction, 
$\sfM \cap \cF^{\ell-1} \simeq \sfM_{\le a}$ is a minimal space.  So it suffices to show that $\sfM$ has the homotopy type of a space obtained from $\sfM \cap \cF^{\ell-1}$ by attaching $b_\ell(\sfM)$ $\ell$-cells, 
where $b_k(\sfM)=\rank H_k(\sfM;\Z)$ denotes the $k$-th Betti number of $\sfM$.

By the Lefschetz hyperplane section theorem of Hamm and L\^{e} \cite{HL} (see also \cite{GM}), $\sfM$ is obtained from $\sfM \cap \cF^{\ell-1}$ by attaching at least $b_\ell(\sfM)$ $\ell$-cells, and the number of $\ell$-cells is equal to the rank of the homology group $H_\ell(\sfM,\sfM\cap\cF^{\ell-1})$.  Using the fact that the Morse function $f$ is perfect repeatedly, we see that the long exact sequence of the pair 
$(\sfM,\sfM\cap\cF^{\ell-1}) \simeq (\sfM_{\le b},\sfM_{\le a})$ splits into short exact sequences as in \eqref{eq:perfect}.  In particular, we have $H_\ell(\sfM) \cong H_\ell(\sfM,\sfM\cap\cF^{\ell-1})$, and $\sfM$ has the homotopy type of a space obtained from the hyperplane section $\sfM\cap\cF^{\ell-1}$ by attaching precisely $b_\ell(\sfM)$ $\ell$-cells.
\end{proof}

A similar proof of minimality was recently given by Yoshinaga \cite{Yo}.

\section{Local Systems}
\label{sec:local systems}
As noted in the introduction, the cohomology of the complement of a complex hyperplane arrangement with coefficients in a (complex) local system is of interest in the study of multivariable hypergeometric integrals, among other applications.

Let $\A$ be a hyperplane arrangement in the complex vector space $V\cong \C^\ell$.
Let $\rho:\pi_1(\sfM) \to \GL_m(\C)$ be a complex representation of the fundamental group of the complement $\sfM$ of $\A$, and denote by $\LL$ the corresponding
rank $m$ local system of coefficients on $\sfM$.  For such a local system, stratified Morse theory was
used in \cite{C1} to construct a complex $(K^\bul(\A),\D^\bul)$,
the cohomology of which is naturally isomorphic to
$H^\bul(\sfM;\LL)$.  We  recall this construction brieflyþ.

Let $\cF$ be a complete flag in $V$ which is transverse to the stratification determined by $\A$ as in \eqref{eq:flag}.  
Let $\sfM^{q} = \cF^q \cap \sfM$ for each $q$.  Let
$K^q=H^q(\sfM^q,\sfM^{q-1};\LL)$, and denote by $\D^q$ the
boundary homomorphism $H^{q}(\sfM^{q},\sfM^{q-1};\LL) \to
H^{q+1}(\sfM^{q+1},\sfM^{q};\LL)$ of the triple
$(\sfM^{q+1},\sfM^q,\sfM^{q-1})$.  The following compiles several
results from \cite{C1}.

\begin{thm}\label{theorem:Kdot}
Let $\LL$ be the complex local system on $\sfM$ corresponding to the representation 
$\rho:\pi_1(\sfM) \to \GL_m(\C)$.

\begin{enumerate}
\item \label{item:Kdot1}
For each $q$, $0\le q \le \ll$, we have
$H^i(\sfM^q,\sfM^{q-1};\LL)=0$ if $i \neq q$, and $\dim_\cn
H^q(\sfM^q,\sfM^{q-1};\LL) = m\cdot b_q(\sfM)$.

\item \label{item:Kdot2}
The system of complex vector spaces and linear maps
$(K^\bul,\D^\bul)$,
\[
K^0 \xrightarrow{\ \D^{0}\ } K^1 \xrightarrow{\ \D^1\ } K^2
\xrightarrow{\phantom{\D^{1}}} \cdots
\xrightarrow{\phantom{\D^{1}}} K^{\ll-1} \xrightarrow{\
\D^{\ll-1}\,} K^\ll,
\]
is a complex $(\D^{q+1}\circ\D^q=0)$.  The cohomology of this
complex is naturally isomorphic to $H^\bul(\sfM;\LL)$, the
cohomology of $\sfM$ with coefficients in $\LL$. \qed
\end{enumerate}
\end{thm}

\begin{cor}[\cite{C2}] \label{cor:morse inequalities}
For the rank $m$ local system $\LL$, let $\beta_q=\dim_\C H^q(\sfM;\LL)$, and write $b_q=b_q(\sfM)$.
Then, for $0\le q \le \ll$, we have
\begin{equation*} \label{eq:weak}
\beta_q \le m\cdot b_q,
\end{equation*}
and
\begin{equation*} \label{eq:strong}
\beta_q - \beta_{q-1} + \dots \pm \beta_0 \le m\cdot(b_q - b_{q-1} + \dots \pm b_0). \qed
\end{equation*}
\end{cor}
These are the weak and strong Morse inequalities arising from the complex $(K^\bul,\Delta^\bul)$ since $\dim_\C K^q = m\cdot b_q$.  In particular, for any complex local system, the cohomology groups $H^q(\sfM;\LL)$ are finite dimensional, resolving a question raised by Aomoto and Kita \cite{AoK1} in the context of rank one local systems.

\begin{rem}
Let $W$ be the minimal CW-complex resulting from (inductive) application of Theorem \ref{thm:min}.
The complex $(K^\bul,\Delta^\bul)$ may be realized as the cellular (co)chain complex of $W$ with coefficients in the local system $\LL$.
\end{rem}

In the rest of this paper we  focus on rank one local systems.  
Let $\bl=(\la_1, \ldots,\la_n)$ be a set of complex weights for
the hyperplanes of $\A$. 
Let $t_j= \exp(-2\pi\ii\la_j)$ and
$\b{t}=(t_1,\dots,t_n)\in (\cn^*)^n$. Associated to $\bl$, we have
a rank one representation $\rho:\pi_1(\sfM) \to \cn^*$, given by
$\rho(\c_j) =t_j$, where $\c_j$ is any meridian loop about the
hyperplane $H_j$ of $\A$, and a corresponding rank one local
system $\LL=\LL_{\b{t}}=\LL_{\bl}$ on $\sfM$.  Note that weights
$\bl$ and $\bl'$ yield identical representations and local systems
if $\bl-\bl'\in\Z^n$.

The dimensions of the terms, $K^{q}$, of the complex
$(K^\bul,\D^\bul)$ are independent of the local system 
$\LL$.  For a rank one local system, they are given by $\dim K^q = b_q(\sfM)$.  In this context, 
write $\D^\bul=\D^\bul(\b{t})$ to indicate the dependence
of the complex on $\b{t}$, and view these boundary maps as
functions of $\b{t}$.  Let $\Lambda=\cn[x_1^{\pm 1},\dots,x_n^{\pm 1}]$
be the ring of complex Laurent polynomials in $n$ commuting
variables.

\begin{thm}[\cite{CO1}] \label{thm:univcx}
For an arrangement $\A$ of $n$ hyperplanes with complement $\sfM$,
there exists a universal complex $(\sK^\bul,\D^\bul(\b{x}))$ with
the following properties:
\begin{enumerate}
\item \label{item:univcx1}
The terms are free $\Lambda$-modules, whose ranks are given by the
Betti numbers of $\sfM$, $\sK^q \simeq \Lambda^{b_q(\A)}$.

\item \label{item:univcx2}
The boundary maps, $\D^q(\b{x}): \sK^q \to \sK^{q+1}$ are
$\L$-linear.

\item \label{item:univcx3}
For each $\b{t}\in (\cn^*)^n$, the specialization $\b{x} \mapsto
\b{t}$ yields the complex $(K^\bul,\D^\bul(\b{t}))$, the
cohomology of which is isomorphic to $H^\bul(\sfM;\LL_\b{t})$, the
cohomology of $\sfM$ with coefficients in the local system
associated to $\b{t}$. \qed
\end{enumerate}
\end{thm}
The entries of the boundary maps $\D^{q}(\b{x})$ are elements of
the Laurent polynomial ring $\L$, the coordinate ring of the
complex algebraic $n$-torus.  Via the specialization $\b{x}
\mapsto \b{t} \in (\cn^*)^n$, we view them as holomorphic
functions $(\cn^*)^n\to\cn$.  Similarly, for each $q$, we view
$\D^{q}(\b{x})$ as a holomorphic map
$\D^{q}:(\cn^*)^n\to\Mat(\cn)$, $\b{t}\mapsto \D^{q}(\b{t})$
from the complex torus to matrices with complex entries.

\begin{rem}
Let $W$ be the minimal CW-complex resulting from application of Theorem \ref{thm:min}, and let $\widetilde{W}$ be the universal cover of $W$.  
The complex $(\sK^\bul,\Delta^\bul(\b{x}))$ may be realized as 
$\Hom^{G}(C_\bul(\widetilde{W}),\Lambda)$, where $G=\pi_1(W)=\pi_1(\sfM)$, $C_\bul(\widetilde{W})$ is the (cellular) chain complex of $\widetilde{W}$, and $\Lambda\cong \C[\Z^n]$ is the $G$-module corresponding to the action of $G$ on the abelianization $G/[G,G] = H_1(W) = \Z^n$ by (left) translation.
In \cite{DP2}, Dimca and Papadima show that the complex $C_\bul(\widetilde{W})$ is itself an invariant of the arrangement $\A$.
\end{rem}

The universal complex $\sK^\bul$ is closely related to another universal complex defined by Aomoto \cite{Aom1} using the Orlik-Solomon algebra 
$A(\A)$.  This graded algebra, isomorphic to the cohomology $H^*(\sfM;\C)$ (see \cite{OrS1,OrT1}), is the quotient of the exterior
algebra $E(\A)$ generated by 1-dimensional
classes $e_j$, $1\leq j \leq n$, by a homogeneous ideal $I(\A)$. 
Let $[n]=\{1,\ldots,n\}$. Refer to the hyperplanes by their subscripts and order
them accordingly. Given $S \subset [n]$, denote the flat $\bigcap_{j \in S} H_j$ by $\cap{S}$.  
If $\cap{S} \neq\emptyset$, call $S$ independent if the codimension of $\cap{S}$ in $V$ is equal to $|S|$, and dependent 
if $\codim(\cap{S}) < |S|$.    If $S=(j_1,j_2,\dots,j_q)$, let 
$e_S = e_{j_1} e_{j_2} \cdots e_{j_q}$ denote the corresponding basis element 
of the exterior algebra.  Define $\partial{e_S} = \sum_{p=1}^q (-1)^{p-1} e_{S \setminus \{j_p\}}$.  The ideal $I(\A)$ is generated by
\[
\{ \partial{e_S} \mid S\text{ is dependent}\} \bigcup
\{ e_S \mid \cap{S} = \emptyset\}. 
\]

For $S \subset [n]$, let 
$a_S$ denote the image of $e_S$ in $A(\A) = E(\A)/I(\A)$.
The algebra $A(\A)$ has a $\cn$-basis called the $\nbc$ basis.
A subset 
$S$ of $[n]$ is a {\em circuit} if it is a minimally dependent set: $S$ is dependent
but every nontrivial subset of $S$ is independent. 
Call $T=(j_1< \cdots <j_p)\subset [n]$ a broken circuit if
there exists $k\in [n]$ so that $k< j_1$ and $(k,T)$ is a circuit. The $\nbc$ basis
consists of all elements $a_S$ of $A(\A)$ corresponding to subsets $S$ of $[n]$ which contain no broken circuit \cite{OrT1}. 
  
Let $a_\bl=\sum_{j=1}^n  \la_j a_j \in A^1(\A)$ and
note that $a_\bl a_\bl=0$ because 
$A(\A)$ is a quotient of an exterior algebra.
Thus we have a complex $(A^{\bul}(\A), a_{\bl})$. Let
$\by=\{y_1,\ldots,y_n\}$ be a set of indeterminates in one-to-one
correspondence with the hyperplanes of $\A$. Let $R=\cn[\by]$ be
the polynomial ring in $\by$. 
Define a graded $R$-algebra:
$\sA^{\bul}=\sA^{\bul}(\A)=R\otimes_{\cn}A^{\bul}(\A)$.
 Let
 $a_{\by}=\sum_{j=1}^n y_j\otimes a_j \in \sA^1$.
The complex $(\sA^{\bul}(\A), a_{\by})$
\begin{equation}
\label{Aomotocomplex} 0 \ra \sA^0(\A) \xrightarrow{a_{\by}}
\sA^1(\A) \xrightarrow{a_{\by}}
 \dots
\xrightarrow{a_{\by}}
  \sA^{\ell}(\A)
\ra 0
\end{equation}
is called the {\em Aomoto complex}. Its 
specialization $\b{y}\mapsto \bl$  is the complex $(A^{\bul}(\A), a_{\bl})$.

\begin{thm}[\cite{CO1}] \label{thm:approx}
For any arrangement $\A$, the Aomoto complex
$(\sA^{\bul},a_{\b{y}})$ is chain equivalent to the linearization
of the universal complex $(\sK^{\bul},\D^{\bul}(\b{x}))$ at the point 
$\bone=(1,\dots,1)  \in (\cn^*)^n$. \qed
\end{thm}

For certain classes of arrangements, the boundary maps of the universal complex $(\sK^{\bul},\D^{\bul}(\b{x}))$ may be described explicitly.  See, for instance, Hattori \cite{Ha} for general position arrangements.  
In the case where 
the arrangement is defined by real equations, progress on this problem has 
been recently made by Yoshinaga \cite{Yo}.  However, for 
an arbitrary arrangement, these boundary maps are not known.  
Consequently, while the complex $(K^{\bul},\D^{\bul}(\b{t}))$ computes local system
cohomology in principle, we do not know how to calculate the
groups $H^q(\sfM;\LL_{\b{t}})$ explicitly for arbitrary weights. 

It is an interesting question to
determine  the stratification of  the
space of all weights with respect to the
local system cohomology groups.  Each point $\b{t}\in\TT$ gives rise to a
local system $\LL=\LL_{\b{t}}$ on the complement $\sfM$.  
Define the {\em characteristic varieties}
\[
\Sigma^{q}_m(\sfM)=\{\b{t}\in\TT \mid \dim
H^{q}(\sfM;\LL_{\b{t}})\ge m\}.
\]
These loci are algebraic subvarieties of $\TT$, which are
invariants of the homotopy type of $\sfM$.  See Arapura \cite{Ar}
and Libgober \cite{Li} for detailed discussions of these varieties
in the contexts of quasiprojective varieties and plane algebraic
curves. The characteristic varieties are closely related to the
resonance varieties.

Each point $\bl\in\cn^{n}$ gives rise to an element
$a_{\bl}\in A^{1}$ of the Orlik-Solomon algebra $A^\bul=A^\bul(\A)$.  Define the
{\em resonance varieties}
\[
\cR^{q}_m(A)=\{\b{\la}\in\cn^{n} \mid \dim
H^{q}(A^{\bul},a_\bl)\ge m\}.
\]
These subvarieties of $\cn^{n}$ are invariants of the
Orlik-Solomon algebra $A(\A)$.  See Falk~\cite{Fa} and
Libgober and Yuzvinsky \cite{LY} for detailed discussions of these
varieties. 

\begin{thm}[\cite{CO1}] \label{thm:tcone}
Let $\A$ be an arrangement in $\cn^\ll$ with complement $\sfM$ and
Orlik-Solomon algebra $A^\bul$.  For each $q$ and $m$, the resonance
variety $\cR_{m}^{q}(A)$ coincides with the tangent cone of the
characteristic variety $\Sigma_{m}^{q}(\sfM)$ at the point
$\bone=(1,\dots,1)\in \TT$. \qed
\end{thm}

The characteristic varieties are known to be unions of
torsion-translated subtori of $\TT$, see \cite{Ar}.  In
particular, all irreducible components of $\Sigma_{m}^{q}(\sfM)$
passing through $\bone$ are subtori of $\TT$.  Consequently, all
irreducible components of the tangent cone are linear subspaces of
$\cn^{n}$.

\begin{cor} \label{cor:falkconj}
For each $q$ and $m$, the resonance variety $\cR^{q}_{m}(A)$ is
the union of an arrangement of subspaces in $\cn^{n}$. \qed
\end{cor}

For $q=1$, these results were established by Cohen and Suciu
\cite{CS4}, see also Libgober and Yuzvinsky \cite{Li,LY}.  For the
discriminantal arrangements of Schechtman and Varchenko \cite{SV},
they were established in \cite{C3}. In particular, as conjectured
by Falk \cite[Conjecture 4.7]{Fa}, the resonance varieties
$\cR_{m}^{q}(A)$ were known to be unions of linear subspaces in
these instances. Corollary \ref{cor:falkconj} above resolves this
conjecture positively for all arrangements in all dimensions.
Theorem \ref{thm:tcone} and Corollary \ref{cor:falkconj} have been
obtained by Libgober in a more general situation, see \cite{L2}.

There are examples of arrangements for which the characteristic
varieties contain (positive dimensional) components which do not
pass through $\bone$ and hence cannot be detected by the resonance
variety, see Suciu \cite{Su}. In some of these cases, the local system 
cohomology is nontrivial, while the cohomology of the Orlik-Solomon 
complex vanishes.

\section{Moduli Spaces}
\label{sec:moduli}

In the rest of the paper we pass from consideration of a fixed 
arrangement to the study of all arrangements of a given combinatorial
type.
Fix a pair $(\ell,n)$ with $n\ge \ell\geq 1$ and consider
families of essential  $\ell$-arrangements with $n$ linearly
ordered hyperplanes. In order to define the notions of
combinatorial type and degeneration, we must allow for
the coincidence of several hyperplanes. We call these new objects
{\em multi-arrangements}.  If there is no coincidence, we call
the arrangement \emph{simple}.

Introduce coordinates $u_1,\ldots,
u_{\ell}$ in $V$ and  choose a degree one polynomial $\al_j
=b_{j,0}+\sum_{k=1}^\ell b_{j,k}u_k$  for the
hyperplane $H_j \in \A$ so $H_j$ is defined by 
$\al_j=0$. Note that $\al_j$ is unique up to a constant. 
Embed $V$ in projective space $ \CP^\ell$ and call the
complement of $V$ the infinite hyperplane, $H_{n+1}$, defined by $u_0=0$.
 We call $\Ai=\A \bigcup
H_{n+1}$ the {\em projective closure} of $\A$. It is
an arrangement in $\CP^\ell$. Give $H_{n+1}$ the weight
$\la_{n+1}=-\sum_{j=1}^n \la_j$.  
We agree that
the hyperplane at infinity, $H_{n+1}$,  is largest in the ordering.
We may therefore view the
projective closure of the arrangement as an $(n+1)\times (\ell+1)$
matrix of complex numbers
\begin{equation}\label{eq:point}
\sfb=
\begin{pmatrix}
b_{1,0} & b_{1,1} & \cdots & b_{1,\ell}\\
b_{2,0} & b_{2,1} & \cdots & b_{2,\ell}\\
\vdots  & \vdots  & \ddots & \vdots \\
b_{n,0} & b_{n,1} & \cdots & b_{n,\ell}\\
1       & 0       & \cdots & 0
\end{pmatrix}
\end{equation}
whose rows correspond to the hyperplanes of $\Ai$. Thus $(\cp^\ell)^{n}$ may
be viewed as the moduli space of all ordered  multi-arrangements
in $\cp^\ell$ with $n$ hyperplanes together with the hyperplane at
infinity.

Given an
arrangement $\A$, the set $S=(j_1,\ldots,j_q)$ is  dependent (in the projective closure) 
if the corresponding row vectors of \eqref{eq:point} are linearly dependent.
Let $\Dep(\A)_q$ be the set of dependent $q$-tuples and let
$\Dep(\A)=\bigcup_{q>1} \Dep(\A)_q$.
Two essential simple
arrangements are  combinatorially equivalent if and only if they
have the same dependent sets. We call $\cT$ their {\em
combinatorial type} and write $\Dep(\cT)$. Note that an arbitrary
collection of subsets of $[n+1]$ is not necessarily  realizable as a
dependent (or independent) set. For example, the collection
$\{123, 124,134\}$ is not realizable as a dependent set, since
these dependencies imply the dependence of 234.

The combinatorial type is, in fact, determined by $\Dep(\cT)_{\ell+1}$, 
see, for instance, Terao \cite{T2}.  
Given a subset $J\subset [n+1]$ of
cardinality $\ell+1$, let $\Delta_J(\sfb)$ denote the determinant
of the $(\ell+1)\times(\ell+1)$ submatrix of $\sfb$ whose rows are specified by
$J$. Given a realizable type $\cT$, the moduli space of type $\cT$
is
\[
\sfX(\cT)=\{\sfb \in (\cp^{\ell})^{n} \mid \Delta_{J}(\sfb) = 0
\text{ for } J\in \Dep(\cT)_{\ell+1},~\Delta_{J}(\sfb)\neq0
\text{ else} \}.
\]
If $\cG$ is the type of a general position arrangement, then
$\Dep(\cG)=\emptyset$
 and the moduli space $\sfX(\cG)$ is a dense, open subset of $(\cp^\ell)^n$. 
Define a partial order on combinatorial
types as follows: $\cT \ge \cT' \iff \Dep(\cT) \subseteq
\Dep(\cT')$. The combinatorial type $\cG$  is the maximal element
with respect to this partial order. Write $\cT > \cT'$ if
$\Dep(\cT) \subsetneq \Dep(\cT')$. If $\cT > \cT'$, we say that
$\cT$ {\em covers} $\cT'$ and $\cT'$ is a {\em degeneration} of
$\cT$ if there is no realizable combinatorial type $\cT''$ with
$\cT>\cT''>\cT'$. In this case we define the relative dependence
set
\[
\Dep(\cT',\cT)= \Dep(\cT') \setminus \Dep(\cT).
\]

 Let
\[
\sfY(\cT) = \{\sfb \in (\cp^{\ell})^{n} \mid \Delta_{J}(\sfb) \neq
0 \text{ for } J \notin \Dep(\cT)_{\ell+1}\}.
\]
Then the moduli space of type $\cT$ may be realized as
\[
\sfX(\cT)=\{\sfb \in \sfY(\cT) \mid \Delta_{J}(\sfb)=0 \text{ for
} J \in \Dep(\cT)_{\ell+1}\}.
\]
Note that $\sfX(\cG)=\sfY(\cG)$.  For any other type $\cT$, the moduli
space $\sfX(\cG)$  may be
realized as
\[
\sfX(\cG)=\{\sfb \in \sfY(\cT) \mid \Delta_{J}(\sfb) \neq 0 \text{
for } J \in \Dep(\cT)_{\ell+1}\}.
\]
If $\cT \neq \cG$, then $\sfX(\cT)$ and $\sfX(\cG)$ are disjoint
subspaces of $\sfY(\cT)$.  Let $i_{\cT}:\sfX(\cG) \to \sfY(\cT)$
and $j_{\cT}:\sfX(\cT) \to \sfY(\cT)$ denote the natural
inclusions. We showed in \cite{CO4} that for any combinatorial
type $\cT$, the inclusion $i_\cT:\sfX(\cG) \to \sfY(\cT)$ induces
a surjection $(i_\cT)_*:H_1(\sfX(\cG)) \to H_1(\sfY(\cT))$.

For the type $\cG$ of general position arrangements, the closure
of the moduli space is $\overline{\sfX}(\cG)= (\cp^\ll)^n$.   The
divisor $\sfD(\cG)=\overline{\sfX}(\cG)\setminus \sfX(\cG)$ is
given by $\sfD(\cG) = \bigcup_J \sfD_J$, whose components, $\sfD_J
=\{\sfb \in (\cp^\ll)^n \mid \D_J(\sfb)=0\}$, are irreducible
hypersurfaces indexed by $J=\{j_1,\dots,j_{\ll+1}\}$.
Choose a basepoint $\b{c} \in \sfX(\cG)$, and for each $\ll+1$
element subset $J$ of $[n+1]$, let $\b{d}_J$ be a generic point in
$\sfD_J$. Let $\Gamma_{\!J}$ be a meridian loop based at $\b{c}$
in $\sfX(\cG)$ about the point $\b{d}_J \in \sfD_J$.  Note that
$\b{c} \in \sfY(\cT)$ and that $\Gamma_{\!J}$ is a (possibly
null-homotopic) loop in $\sfY(\cT)$ for any combinatorial type
$\cT$. We showed in \cite{CO4} that for any combinatorial type
$\cT$, the homology group $H_1(\sfY(\cT))$ is generated by the
classes $\{[\Gamma_{\!J}] \mid J \not\in \Dep(\cT)_{\ell+1}\}$.
In particular, the homology group $H_1(\sfX(\cG))$ is generated by the
classes $[\Gamma_{\!J}]$, where $J$ ranges over all $\ll+1$
element subsets of $[n+1]$.

It is easy to see that  the moduli space $\sfX(\cT')$ has complex
codimension one in the closure $\overline{\sfX}(\cT)$  if and only
if $\cT$ covers $\cT'$. The next theorem is essential for later
results.

\begin{thm}[\cite{CO4}] \label{thm:lincomb}
\label{thm:comp} Let $\cT$ be a combinatorial type which covers
the type $\cT'$. Let $\sfb'$ be a point in $\sfX(\cT')$, and $\c
\in \pi_1(\sfX(\cT),\sfb)$ a simple loop in $\sfX(\cT)$ about
$\sfb'$. Then the homology class $[\c]$ satisfies
\begin{equation} \label{eq:compat}
(j_\cT)_*([\c]) = \sum_{J \in \Dep(\cT',\cT)} m_J \cdot
[\Gamma_{\!J}],
\end{equation}
where $m_J$ is the order of vanishing of the restriction of
$\Delta_J$ to $\overline{\sfX}(\cT)$ along $\sfX(\cT')$. \qed
\end{thm}

\section{Gauss-Manin Connections}
\label{sec:GM}

The moduli space  $\sfX(\cT)$ is not necessarily connected. The
existence of a combinatorial type whose moduli space has at least
two components follows from examples of Rybnikov \cite{Ry}.
Let $\sfB(\cT)$ be a smooth component of the moduli space.
Corresponding to each $\sfb\in \sfB(\cT)$, we have an arrangement
$\A_\sfb$, combinatorially equivalent to $\A$, with hyperplanes
defined by the first $n$ rows of the matrix equation $\sfb \cdot
\tilde{\sf u}=0$, where $\tilde{\sf u} =
\begin{pmatrix} 1 & u_1 & \cdots & u_\ell\end{pmatrix}^\top$.
Let $\sfM_\sfb=M(\A_\sfb)$ be the complement of $\A_\sfb$. Let
\[
\sfM(\cT) = \{ (\sfb,{\sf u}) \in (\cp^\ell)^n \times \cn^\ell
\mid \sfb \in \sfB(\cT) \ \hbox{and}\ {\sf u} \in \sfM_\sfb\},
\]
and define $\pi_{\cT}:\sfM(\cT) \to \sfB(\cT)$ by
$\pi_{\cT}(\sfb,{\sf u})=\sfb$.  Since $\sfB(\cT)$ is connected by
assumption, a result of Randell \cite{Ra1} implies that
$\pi_{\cT}:\sfM(\cT) \to \sfB(\cT)$ is a bundle, with fiber
$\pi_{\cT}^{-1}(\sfb) = \sfM_\sfb$.

For each $\sfb \in\sfB(\cT)$, weights $\bl$ define a local system
$\LL_\sfb$ on $\sfM_\sfb$.  Since $\pi_{\cT}:\sfM(\cT) \to
\sfB(\cT)$ is locally trivial, there is an associated flat vector
bundle $\pi^{q}:\b{H}^{q}(\LL)\to\sfB(\cT)$, with fiber
$(\pi^{q})_\LL^{-1}(\sfb)=H^{q}(\sfM_{\sfb};\LL_{\sfb})$ at
$\sfb\in\sfB(\cT)$ for each $q$, $0\le q \le \ell$.   Fixing a
basepoint $\sfb\in\sfB(\cT)$, the operation of parallel
translation of fibers over curves in $\sfB(\cT)$  provides a
complex representation
\begin{equation} \label{eq:Hqrep}
\Psi_{\cT}^{q}:\pi_{1}(\sfB(\cT),\sfb) \longrightarrow
\Aut_\cn(H^{q}(\sfM_{\sfb};\LL_{\sfb})).
\end{equation}
The cohomology of the Morse theoretic complex $K^\bul(\A_\sfb)$ 
is isomorphic to the cohomology of $\sfM_\sfb$ with coefficients in 
the local system $\LL_{\sfb}$.  
The fundmental
group of $\sfB(\cT)$ acts by chain automorphisms on this complex,
see \cite[Cor.~3.2]{CO3}, yielding a representation
\begin{equation*} \label{eq:KArep}
\psi_{\cT}^{\bul}:\pi_{1}(\sfB(\cT),\sfb) \longrightarrow
\Aut_\cn(K^{\bul}(\A_{\sfb})).
\end{equation*}

\begin{thm} \label{thm:inducedrep}
The representation $\Psi^{q}_\cT:\pi_{1}(\sfB(\cT),\sfb) \to
\Aut_\cn(H^{q}(\sfM_{\sfb};\LL_{\sfb}))$ is induced by the
representation $\psi^{\bul}_\cT:\pi_{1}(\sfB(\cT),\sfb) \to
\Aut_\cn(K^{\bul}(\A_{\sfb}))$. \qed
\end{thm}

The vector bundle $\pi^q:\b{H}^q(\LL)\to\sfB(\cT)$ supports a
Gauss-Manin connection corresponding to the representation
(\ref{eq:Hqrep}).  Over a manifold $X$, there is a well known
equivalence between complex local systems and complex vector bundles
equipped with flat connections, see \cite{De1,Ko}.  Let $\b{V}\to
X$ be such a bundle, with connection $\nabla$.  The latter is a
$\cn$-linear map $\nabla:\cE^0(\b{V}) \to \cE^1(\b{V})$, where
$\cE^p(\b{V})$ denotes the complex $p$-forms on $X$ with values in
$\b{V}$, which satisfies $\nabla(f\sigma)= \sigma df + f
\nabla(\sigma)$ for a function $f$ and $\sigma\in\cE^0(\b{V})$.
The connection extends to a map $\nabla:\cE^p(\b{V}) \to
\cE^{p+1}(\b{V})$ for $p\ge 0$, and is flat if the curvature
$\nabla\circ\nabla$ vanishes.  Call two connections $\nabla$ and
$\nabla'$ on $\b{V}$ isomorphic if $\nabla'$ is obtained from
$\nabla$ by a gauge transformation, $\nabla'=g\circ\nabla\circ
g^{-1}$ for some $g:X\to\Hom(\b{V},\b{V})$.

The aforementioned equivalence is given by $(\b{V},\nabla) \mapsto
\b{V}^{\nabla}$, where $\b{V}^{\nabla}$ is the local system, or
locally constant sheaf, of horizontal sections $\{\sigma \in
\cE^0(\b{V})\mid \nabla(\sigma)=0\}$.  There is also a well known
equivalence between local systems on $X$ and finite dimensional
representations of the fundamental group of $X$.  Note
that isomorphic connections give rise to the same representation.
Under these equivalences, the local system on $X=\sfB(\cT)$
induced by the representation $\Psi^q_\cT$ corresponds to a flat
connection on the vector bundle $\pi^q:\b{H}^q(\LL)\to\sfB(\cT)$,
the Gauss-Manin connection.

Let $\c\in\pi_1(\sfB(\cT),\sfb)$, and let $g:\bS^1 \to \sfB(\cT)$
be a representative loop.  Pulling back the bundle
$\pi^q:\b{H}^q(\LL) \to \sfB(\cT)$ and the Gauss-Manin connection
$\nabla$, we obtain a flat connection $g^*(\nabla)$ on the vector
bundle over the circle corresponding to the representation of
$\pi_1(\bS^1,1)=\langle\zeta\rangle=\Z$ given by $\zeta \mapsto
\Psi^q_\cT(\c)$.  This vector bundle is trivial since any map from
the circle to the relevant classifying space is null-homotopic.
Specifying the flat connection $g^*(\nabla)$ amounts to choosing a
logarithm of $\Psi^q_\cT(\c)$.  The connection $g^*(\nabla)$ is
determined by a connection $1$-form $dz/z \otimes
\Omega^q_\cT(\c)$, where the connection matrix $\Omega^q_\cT(\c)$
corresponding to $\c$ satisfies $\Psi^q_\cT(\c) = \exp(-2 \pi\ii
\Omega^q_\cT(\c))$.  If $\c$ and $\hat\c$ are conjugate in
$\pi_1(\sfB(\cT),\sfb)$, then the resulting connection matrices
are conjugate, and the corresponding connections on the trivial
vector bundle over the circle are isomorphic.  In this sense, the
connection matrix $\Omega^q_\cT(\c)$ is determined by the homology
class $[\c]$ of $\c$.

In the special case when $\cT$ covers $\cT'$ and
$\c\in\pi_1(\sfB(\cT),\sfb)$ is a simple loop linking $\sfB(\cT')$
in $\overline{\sfB(\cT)}$, we denote the corresponding Gauss-Manin
connection matrix in the bundle $\pi^q:\b{H}^q(\LL) \to \sfB(\cT)$ by
$\Omega^q_{\LL}(\sfB(\cT'),\sfB(\cT))$.  
The relationship between the homology classes of the loop $\c$ and 
the loops $\Gamma_J$ in the moduli space of a general position 
arrangement exhibited in Theorem \ref{thm:lincomb} suggests 
an analogous relationship between the corresponding 
Gauss-Manin endomorphisms.  

For nonresonant weights $\bl$, 
this relationship is pursued in \cite{CO4}.  In this situation, the local system cohomology 
is concentrated in the top dimension, and is isomorphic to the cohomology of the Orlik-Solomon complex, $H^\ell(\sfM;\LL_\bl) \cong H^\ell(A^\bul(\A),a_\bl)$.  
Moreover, there is a surjection $P:H^\ell(A(\cG),e_\bl) \twoheadrightarrow 
H^\ell(A^\bul(\A),a_\bl)$ from the cohomology of the Orlik-Solomon complex of a general position arrangement to that of $\A$, see \cite[Theorem 6.5]{CO4}.

\begin{thm} \label{thm:nonresGM}
Let $\cT$ be a combinatorial type which covers the type $\cT'$.  Let
$\bl$ be a collection of weights which are nonresonant for type $\cT$ 
(and hence for type~$\cG$).  Then the
Gauss-Manin endomorphism $\Omega^\ell_{\LL}(\sfB(\cT'),\sfB(\cT))$ is determined by the
equation
\[
P \cdot \Omega^\ell_{\LL}(\sfB(\cT'),\sfB(\cT)) = \Bigl( \sum_{J \in \Dep_{\ell+1}(\cT',\cT)}
m_J \cdot \Omega^\ell_{\LL}(\sfB(\cT_J),\sfB(\cG))\Bigr) \cdot P,
\]
where $\cT_J$ is the combinatorial type of an arrangement for which $J$ is 
the only dependent set of size $\ell+1$, and $\Omega^\ell_{\LL}(\sfB(\cT_J),\sfB(\cG))
\in \End H^\ell(A(\cG),e_\bl)$ is the corresponding Aomoto-Kita Gauss-Manin connection matrix. \qed
\end{thm}

For arbitrary weights $\bl$, Theorems \ref{thm:lincomb} and \ref{thm:nonresGM} 
motivated the construction of formal connections in 
the Aomoto complex of a general position arrangement in \cite{CO5}.  
These are discussed in Section \ref{sec:formal}.

The Gauss-Manin  connection in local system cohomology has
combinatorial analogs. We have the vector bundle  $\b{A}^q \to
\sfB(\cT)$, whose fiber at $\sfb$ is $A^q(\A_\sfb)$, the $q$-th
graded component of the Orlik-Solomon algebra of the arrangement
$\A_\sfb$. The \textbf{nbc} basis provides a global trivialization
of this bundle. Given weights $\bl$, the cohomology of the complex
$(A^{\bul}(\A_\sfb),a_\bl)$ gives rise to  the flat vector bundle
$\b{H}^q(A) \to \sfB(\cT)$ whose fiber at $\sfb$ is the $q$-th
cohomology group of the Orlik-Solomon algebra,
$H^q(A^{\bul}(\A_\sfb),a_\bl)$. Like their topological
counterparts, these algebraic vector bundles admit flat
connections.  If $\cT$ covers $\cT'$, denote the corresponding 
connection matrix in this cohomology bundle by 
$\Omega^q_{A}(\sfB(\cT'),\sfB(\cT))$.

\section{Formal Connections}
\label{sec:formal}

To determine the endomorphisms
$\Omega^q_{\LL}(\sfB(\cT'),\sfB(\cT))$ and $\Omega^q_{A}(\sfB(\cT'),\sfB(\cT))$, 
we define formal
connections in the Aomoto complex, $(\sA^{\bul}(\cG),a_{\b{y}}))$, of the general position
arrangement of $n$
ordered hyperplanes in $\cn^\ell$.  We
embed the arrangement in projective space as described above
and call the resulting type $\cG_\infty$.
The symmetric group $\Sigma_{n+1}$ 
on $n+1$ letters acts on $\sA^\bul(\cG)$, the rank $\ell$
truncation of the exterior algebra
in $n$ variables,  by permuting the hyperplanes 
of $\cG_\infty$, and on $R$ by permuting the variables $y_j$, where 
$y_{n+1}=-\sum_{j=1}^n y_j$.
In the basis $\{e_j \mid 1\le j \le n\}$ for the exterior algebra, the 
action
of $\sigma\in \Sigma_{n+1}$  is given 
by $\sigma(e_i)=e_{\sigma(i)}$ if $\sigma(n+1)=n+1$, and by 
\[
\sigma(e_i)=\begin{cases}
-e_{\sigma(n+1)}&\text{if $\sigma(i)=n+1$,}\\
e_{\sigma(i)} - e_{\sigma(n+1)}&\text{if $\sigma(i) \neq n+1$,}
\end{cases}
\]
if $\sigma(n+1)\neq n+1$.  
Denote the induced action on the Aomoto complex by 
$\phi_\sigma:\sA^\bul(\cG) \to \sA^\bul(\cG)$, 
\[
\phi_\sigma(e_{i_1} \cdots e_{i_p} \otimes f(y_1,\dots,y_n))=
\sigma(e_{i_1})\cdots\sigma(e_{i_p}) \otimes 
f(y_{\sigma(1)},\dots,y_{\sigma(n)}).
\]

\begin{lem} \label{lem:sym action}
For each $\sigma \in \Sigma_{n+1}$, the map $\phi_\sigma$ is a cochain 
automorphism of the Aomoto complex $(\sA^\bul(\cG),e_{\b{y}})$. \qed
\end{lem}

If $T=(i_1,\ldots,i_p)\subset [n]$ is a $p$-tuple, 
then we write $e_T=e_{i_1}\cdots e_{i_p}$. 
Recall that $\partial{e}_T = \sum_{j=1}^p (-1)^{j-1} e_{T \setminus\{i_j\}}$.  
For $j \in [n]$, let 
$(j,T)=(j,i_1,\ldots,i_p)$ be the $(p+1)$-tuple which adds $j$
to $T$ as its first entry. 
For $S=(s_1,\dots,s_k) \subset [n+1]$, let 
$\sigma_S$ denote 
the permutation $\big(\begin{smallmatrix}1 & 2 & \cdots & k\\
s_1 & s_2 & \cdots & s_{k}\end{smallmatrix}\big)$.  
Write $S \equiv T$ if $S$ 
and $T$ are equal sets.  

\begin{definition} \label{def:omega}
Let $T\subset [n]$ be a $p$-tuple, 
 $S \subset [n+1]$ have size  $q+1$, and $j\in[n]$. 
If $S=S_0=[q+1]$, define the endomorphism
$\tilde{\omega}^\bullet_{S_0}:(\sA^{\bul}(\cG),e_{\b{y}}) \rightarrow
(\sA^{\bul}(\cG),e_{\b{y}})$ by
\[
\tilde{\omega}_{S_0}^p(e_T)=
\begin{cases}
y_j \p e_{(j,T)} & \text{if $p=q$ and $S_0\equiv (j,T)$,}\\
e_{\b{y}} \p e_T & \text{if $p=q+1$ and $S_0\equiv T$,}\\
0 & \text{otherwise.}
\end{cases}
\]
If $S\neq S_0$, 
define $\tilde{\omega}^\bullet_S = 
\phi_{\sigma_S}^{} \circ \tilde{\omega}_{S_0}^{\bullet} 
\circ 
\phi_{\sigma_S}^{-1}$.
\end{definition}

\begin{prop}[\cite{CO5}] \label{prop:chain map}
For every subset $S$ of $[n+1]$, the map $\tilde{\omega}^\bullet_S$ is a cochain
homomorphism of the Aomoto complex
$(\sA^{\bul}(\cG),e_{\b{y}}))$. \qed
\end{prop}

The formal connection endomorphisms $\tilde{\omega}^\bullet_S$ are defined 
for the Aomoto complex of the general position type $\cG$. Our aim is to show 
that certain linear combinations of these 
induce endomorphisms of the Aomoto complex of type $\cT$ 
for all pairs of types $\cT'$, $\cT$ 
where $\cT$ covers $\cT'$. This involves multiplicities.
Given $S \subset [n+1]$, let $N_S(\cT)=N_S(\sfb)$ denote the
submatrix of (\ref{eq:point}) with rows specified by $S$. Let
$\rank N_S(\cT)$ be the size of the largest minor with nonzero
determinant. Define the {\em multiplicity} of $S$ in $\cT$ by
\[ 
m_S(\cT)=|S|-\rank N_S(\cT).
\]
It is not hard to see that this definition of multiplicity agrees with the
analytic definition in Theorem \ref{thm:comp}. Let
\[
\tilde{\omega}(\cT',\cT)=\sum_{S\in
\Dep(\cT',\cT)}m_S(\cT')\cdot\tilde{\omega}_S.
\]
For an arrangement $\A$, the Orlik-Solomon algebra 
depends only on the combinatorial type $\cT$, so we write $A(\A)=A(\cT)$.  
If $\A \subset V\cong \C^\ell$, then $A^q(\cT)=0$ for $q > \ell$.  
It follows that $A(\cT)$ may be realized as a quotient of the rank $\ell$ 
truncation of the exterior algebra $E(\A)$, which is itself the Orlik-Solomon algebra 
$A(\cG)$ of the combinatorial type of a general position arrangement.
Denote the rank $\ell$ truncation of the Orlik-Solomon ideal $I(\A)$ by 
$I(\cT) = I(\A) \cap A(\cG)$.
Thus $A(\cT)=A(\cG)/I(\cT)$. The ideal $I(\cT)$ gives rise to a
subcomplex $\sfI^\bul(\cT)$ of the Aomoto complex
$\sA^{\bul}(\cG)$, and we have 
an exact sequence of cochain complexes 
\[
0 \to \sfI^\bul(\cT) \to \sA^{\bul}(\cG) \to \sA^{\bul}(\cT) \to 0.
\]

\begin{thm}[\cite{CO5}]
\label{thm:ideal} If $\cT$ covers $\cT'$, then
$\tilde{\omega}(\cT',\cT)(\sfI^\bul(\cT))\subset \sfI^\bul(\cT)$ so
there is a commutative diagram
\[
\begin{CD}
(\sfI^\bul(\cT),e_{\b{y}}) @>\iota>> (\sA^{\bul}(\cG),e_{\b{y}})
@>p>> (\sA^{\bul}(\cT),a_{\b{y}}) \\
@VV\left.\tilde{\omega}(\cT',\cT)\right|_{\sfI^\bul(\cT)}V
@VV\tilde{\omega}(\cT',\cT)V
@VV{\omega}(\cT',\cT)V \\
(\sfI^\bul(\cT),e_{\b{y}}) @>\iota>> (\sA^{\bul}(\cG),e_{\b{y}})
@>p>> (\sA^{\bul}(\cT),a_{\b{y}}) \\
\end{CD}
\]
where $\iota:\sfI^\bul(\cT) \to \sA^\bul(\cG)$ is the inclusion,
$p:\sA^\bul(\cG) \to \sA^\bul(\cT) = \sA^\bul(\cG)/\sfI^\bul(\cT)$ is the 
natural projection, 
and $\omega(\cT',\cT):\sA^\bul(\cT) \to \sA^\bul(\cT)$ is the induced map. \qed
\end{thm}
We call the map $\omega(\cT',\cT)$ the universal Gauss-Manin endomorphism.

It follows  that for given weights $\bl$, the
specialization $\b{y}\mapsto \bl$ in the chain endomorphism
$\omega(\cT',\cT)$ defines a chain endomorphism
$\omega^\bul_{\bl}(\cT',\cT):A^\bul(\cT)\rightarrow A^\bul(\cT)$. 
Let $\kappa^q = \ker[a_\bl:A^q(\cT) \to A^{q+1}(\cT)]$, and 
write $A^q(\cT) = \kappa^q \oplus A^q(\cT)/\kappa^q$.  
Define $\rho^q:A^q(\cT)\twoheadrightarrow H^q(A^\bul(\cT),a_{\bl})$ 
to be the natural projection $\kappa^q \twoheadrightarrow 
H^q(A^\bul(\cT),a_{\bl})$ on $\kappa^q$, and trivial on $A^q(\cT)/\kappa^q$.  
The map 
$\omega^q_{\bl}(\cT',\cT)$ induces an endomorphism
\[
\Omega^q_{\cC}(\cT',\cT):H^q(A^\bul(\cT),a_{\bl})\to
H^q(A^\bul(\cT),a_{\bl})
\]
determined by the equation $ \rho^q\circ \omega^q_{\bl}(\cT',\cT)
= \Omega^q_{\cC}(\cT',\cT) \circ \rho^q$.

\begin{thm}[\cite{CO5}]
\label{thm:GM1}
Let $\sfM$ be the complement of an arrangement of
type $\cT$ and let $\LL$ be the local system on $\sfM$ defined by
weights $\bl$. Suppose $\cT$ covers $\cT'$.
Then the connection endomorphism
$\Omega^q_{A}(\sfB(\cT'),\sfB(\cT))$ is determined by the equation
\[
\rho^q\circ \omega^q_{\bl}(\cT',\cT) =
\Omega^q_{A}(\sfB(\cT'),\sfB(\cT)) \circ \rho^q.
\]
and hence
$\Omega^q_{A}(\sfB(\cT'),\sfB(\cT))=\Omega^q_{\cC}(\cT',\cT)$. \qed
\end{thm}

Now consider the endomorphisms $\Omega^q_{\LL}(\sfB(\cT'),\sfB(\cT))$ 
of the local system cohomology groups $H^q(\sfM;\LL)$.  
Recall from Theorem~\ref{theorem:Kdot} 
that this cohomology is naturally isomorphic to the cohomology 
of the Morse theoretic complex $(K^\bul(\A),\Delta^\bul)$.  As above, 
let $\varkappa^q=\ker[\Delta^q:K^q(\A) \to K^{q+1}(\A)]$, and 
write $K^q(\A) = \varkappa^q \oplus K^q(\A)/\varkappa^q$.  
Define $\varphi^q:K^q(\A)\twoheadrightarrow H^q(\sfM;\LL)$ 
to be the natural projection $\varkappa^q \twoheadrightarrow 
H^q(\sfM;\LL)$ on $\varkappa^q$, and trivial on $K^q(\A)/\varkappa^q$.  

\begin{thm}[\cite{CO5}]
\label{thm:GM2}
Let $\sfM$ be the complement of an arrangement of
type $\cT$ and let $\LL$ be the local system on $\sfM$ defined by
weights $\bl$. Suppose $\cT$ covers $\cT'$.
Then there is an isomorphism $\tau^q:A^q(\cT)\rightarrow
K^q(\A)$ so that the Gauss-Manin endomorphism
$\Omega^q_{\LL}(\sfB(\cT'),\sfB(\cT))$ in local system cohomology
is determined by the equation
\[
\varphi^q\circ\tau^q \circ \omega^q_{\bl}(\cT',\cT) =
\Omega^q_{\LL}(\sfB(\cT'),\sfB(\cT)) \circ \varphi^q\circ\tau^q. \qed
\]
\end{thm}

\section{Spectrum}
\label{sec:eigen}

The eigenvalues of the Gauss-Manin connection satisfy:
\begin{thm}[\cite{CO3}]
 \label{thm:main1}
The eigenvalues of the universal Gauss-Manin endomorphism 
$\omega^q(\cT',\cT)$ are integral
linear forms in the variables $y_1,\dots,y_n$. Thus
for any system of weights $\bl$, the
eigenvalues of the Gauss-Manin endomorphism in local system cohomology,
$\Omega^q_{\LL}(\sfB(\cT'),\sfB(\cT))$,
are  integral linear combinations of the weights~$\bl$. \qed
\end{thm}

In \cite{CO6} we determined the spectra of these Gauss-Manin endomorphisms.  
Recall the collection $\Dep(\cT)$.
Here it suffices to work with a smaller collection of dependent sets
\[
\Dep(\cT)_q^*=
\{S\in  \Dep(\cT)_q \mid \bigcap_{j \in S}H_j \neq \emptyset\}.
\] 
Let $\Dep(\cT)^*= \bigcup_{q} \Dep(\cT)_q^*$.  If $S \in  \Dep(\cT)^*$,  then
$\codim(\bigcap_{j \in S}H_j) < |S|$. 
If $\cT'$ is a combinatorial type 
for which $\Dep(\cT)^* \subset \Dep(\cT')^*$, let
$\Dep(\cT',\cT)^* = \Dep(\cT')^* \setminus \Dep(\cT)^*$.
If $|S|\geq \ell +2$, then $S \in\Dep(\cT)$ but  $S \in\Dep(\cT)^*$
if and only if every subset of $S$ of cardinality $\ell+1$ is dependent.

Denote the cardinality of $S$ by $s=|S|$.  For 
$1 \le r \le \min(\ell,s-1)$, 
consider the combinatorial type $\cT(S,r)$ defined by
\[
T \in \Dep(\cT(S,r))^* \iff |T \cap S| \geq r+1.
\]
This type is realized by a pencil of hyperplanes indexed by $S$ with a common
subspace of codimension $r$, together with $n-s$ hyperplanes in general 
position.
Note that for $r=1$ the hyperplanes in $S$ coincide, so $\cT(S,r)$   
is a multi-arrangement.

\begin{thm}[\cite{CO6}]
\label{thm:principal}
Let $\cT'$ be a degeneration of a realizable combinatorial type $\cT$.
For each set $S_i \in \Dep(\cT',\cT)^*$, let $r_i$ be minimal 
so that $\Dep(\cT(S_i,r_i))^* \subset  \Dep(\cT')^*$. 
Given the collection 
$\{(S_i, r_i)\}$ there is a unique pair $(S,r)$ with $r=\min \{r_i\}$, 
 $\Dep(\cT(S,r))^* \subset  \Dep(\cT')^*$, and 
 for every  pair $(S_i,r_i)$ where $r_i=r$,  
$S_i \subset S$. \qed
\end{thm}
Let $\cT'$ be a degeneration of $\cT$.
We call the  pair $(S,r)$ 
which satisfies the conditions of Theorem
\ref{thm:principal} the {\em principal  dependence}
of the degeneration. 
Define 
\begin{equation*} \label{eq:(S,r)}
\tilde\omega^\bullet(S,r) = \sum_{K\in\Dep(\cT(S,r))^*} m_K(S,r)\cdot 
\tilde\omega^\bullet_K,
\end{equation*}
where $m_K(S,r)$ is the multiplicity of $K$ in type $\cT(S,r)$.
We showed in the proof of \cite[Thm. 5.1]{CO6}
that the endomorphisms  $\tilde\omega^\bullet(S,r) $ and $\tilde\omega^\bullet(\cT',\cT)$
of $\sA^\bul(\cG)$ induce the same endomorphism in $\sA^\bul(\cT)$,
$\omega^\bullet(S,r) =\omega^\bullet(\cT',\cT)$. It
 follows from
Theorem~\ref{thm:GM2} that for all weights $\bl$, the endomorphism  
$\omega^\bullet_\bl(S,r)$ induces the Gauss-Manin endomorphism 
$\Omega^\bul_{\LL}(\sfB(\cT'),\sfB(\cT))$. Write $\la_S=\sum_{j \in S}\la_j$.

\begin{thm}[\cite{CO6}] \label{thm:diag}
Suppose $\cT$ covers $\cT'$ with principal dependence $(S,r)$. 
Let $\bl$ be a collection of weights satisfying $\la_S \neq 0$.  
Then $\tilde{\omega}^q_\bl(S,r):A^q(\cG) \to A^q(\cG)$, the 
specialization of $\tilde{\omega}^q(S,r)$ at 
$\bl$, is diagonalizable, with eigenvalues $0$ and $\la_S$.  
\begin{enumerate}
\item The $0$-eigenspace has dimension 
\[
\sum_{p=0}^r \binom{s}{p}\binom{n-s}{q-p}-\binom{s-1}{r}\binom{n-s}{q-r}.
\]
\item The $\la_S$-eigenspace has dimension
\[
\sum_{p=r+1}^{\min(q,s)} 
\binom{s}{p}\binom{n-s}{q-p}+\binom{s-1}{r}\binom{n-s}{q-r}. \qed
\]
\end{enumerate}
\end{thm}

Our last
result was stated in \cite{CO6} only for nonresonant
weights but applies in full generality:  
\begin{thm}[\cite{CO6}]
\label{thm:evalues} 
Suppose $\cT$ covers $\cT'$ with principal dependence $(S,r)$. 
Let $\bl$ be a collection of weights satisfying $\la_S \neq 0$.  
Then the Gauss-Manin endomorphism  $\Omega^q_{\LL}(\sfB(\cT'),\sfB(\cT))$
is diagonalizable, with spectrum contained in $\{0,\la_S\}$. \qed
\end{thm}

\section{A Selberg arrangement} \label{sec:selberg}
\begin{figure}[h]
\setlength{\unitlength}{.45pt}
\begin{picture}(300,130)(-200,-110)
\put(-200,-100){\line(0,1){110}}\put(-270,-100){\line(1,1){110}}
\put(-230,-100){\line(0,1){110}}\put(-270,-30){\line(1,0){110}}
\put(-270,-60){\line(1,0){110}}
\put(-155,15){5}\put(-155,-25){4}\put(-155,-55){3}\put(-235,15){1}
\put(-205,15){2}
\put(-224,-122){${\mathcal A}$}

\put(100,-100){\line(0,1){110}}
\put(70,-100){\line(0,1){110}}
\put(30,-60){\line(1,0){110}}
\put(140,-55){345}\put(65,15){1}\put(95,15){2}
\put(75,-122){${\mathcal A}'$}
\end{picture}
\caption{A Selberg Arrangement and One Degeneration}
\label{fig:selberg}
\end{figure}
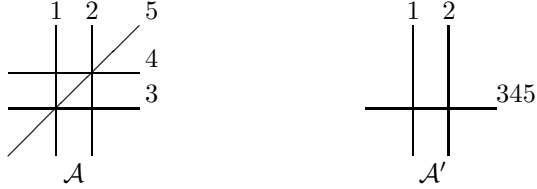

Let $\cS$ be the
combinatorial type of the Selberg arrangement $\A$ in $\C^2$ with
defining polynomial $Q(\A)=u_1u_2(u_1-1)(u_2-1)(u_1-u_2)$
depicted in Figure \ref{fig:selberg}.  
See \cite{Ao,SV,JK} for detailed studies of the Gauss-Manin
connections arising in the context of Selberg arrangements.

Let $\LL$ be the complex rank one local system on the complement $\sfM$ of $\A$ 
corresponding to the point $\b{t}=(t_1,\dots,t_5) \in (\C^*)^5$.  For 
any such local system on the complement of this arrangement, 
there is a choice of weights $\bl=(\la_1,\dots,\la_5) \in \C^5$ so that $t_j=\exp(-2\pi\ii\la_j)$ 
for each $j$, and the local system cohomology $H^*(\sfM;\LL)$ is isomorphic to the 
cohomology of the Orlik-Solomon complex $(A^\bul(\cS),a_\bl)$.  Consequently, 
if $\cS'$ is a degeneration of $\cS$, it suffices to compute the Gauss-Manin endomorphism 
$\Omega^q_A(\sfB(\cS'),\sfB(\cS))=\Omega^q_\LL(\sfB(\cS'),\sfB(\cS))=\Omega^q_\cC(\cS'\cS)$.

Let $\cG$ be the combinatorial type of a general position arrangement of 
five lines in $\C^2$.  
The \textbf{nbc} bases for the Orlik-Solomon algebras $A(\cG)$ and $A(\cS)$ 
give rise to bases for the corresponding Aomoto complexes.  
The Aomoto complex $(\sA^\bul(\cG),e_\b{y})$ is (dual to) the rank two truncation of 
the standard Koszul complex of $\b{y}=y_1,\dots,y_5$ in  the polynomial ring $R=\C[\b{y}]$.  
The Aomoto complex $(\sA^\bul(\cS),a_\b{y})$ of the Selberg arrangement is given by
\[
\sA^0(\cS) \xrightarrow{a_\b{y}} \sA^1(\cS) \xrightarrow{a_\b{y}} \sA^2(\cS),
\]
where $\sA^0(\cG)=R$, $\sA^1(\cG)=R^5$, and $\sA^2(\cG)=R^{6}$.  
Recall that $y_J = \sum_{j \in J} y_j$.  
The boundary maps of this complex have matrices 
\[
\left[\begin{matrix} y_1 & y_2 & y_3 & y_4 & y_5\end{matrix}\right] 
\quad\text{and}\quad
\left[
\begin{matrix}
-y_3 & -y_4 & -y_5 & 0 & 0 & 0 \\
0 & 0 & 0 & -y_3 & -y_4 & -y_5 \\
y_{15} & 0 & -y_5 & y_2 & 0 & 0 \\
0 & y_1 & 0 & 0 & y_{25} & -y_5 \\
-y_3 & 0 & y_{13} & 0 & -y_4 & y_{24}
\end{matrix}
\right].
\]
The projection $p:\sA^\bul(\cG) \to \sA^\bul(\cS)$ is given, in the \textbf{nbc} bases, by
\[
p(e_J)=
\begin{cases}
0 & \text{if $J=\{1,2\}$ or $J=\{3,4\}$,} \\
a_{1,5} - a_{1,3} & \text{if $J=\{3,5\}$,} \\
a_{2,5} - a_{2,4} & \text{if $J=\{4,5\}$,} \\
a_J & \text{otherwise.}
\end{cases}
\]

Let $\cS'$ denote the combinatorial type of the (multi)-arrangement $\A'$ shown in 
Figure \ref{fig:selberg}, a codimension one degeneration of $\cS$.  
The principal dependence of this degeneration is 
$(S,r)$, where $S=345$ and $r=1$.  The corresponding endomorphism 
$\tilde{\omega}^\bul(S,r):\sA^\bul(\cG)\to \sA^\bul(\cG)$ is given by 
\[
\begin{aligned}
\tilde\omega(S,r)&=
\tilde\omega_{34}+\tilde\omega_{35}+\tilde\omega_{45}+
\tilde\omega_{134}+\tilde\omega_{234}+\tilde\omega_{346}+\tilde\omega_{135}+
\tilde\omega_{235}+\tilde\omega_{356}\\
&\qquad+\tilde\omega_{145}+\tilde\omega_{245}
+\tilde\omega_{456}+2\tilde\omega_{345}.
\end{aligned}\]
The matrices of this chain endomorphism are $\tilde\omega^0(S,r)=0$,
\[
\tilde\omega^1(S,r)=
\left[
\begin{matrix}
0 & 0 & 0 & 0 & 0 \\
0 & 0 & 0 & 0 & 0 \\
0 & 0 & y_{45} & -y_4 & -y_5 \\
0 & 0 & -y_3 & y_{35} & -y_5 \\
0 & 0 & -y_3 & -y_4 & y_{34}
\end{matrix}
\right],
\]
\[
\tilde\omega^2(S,r)=
\left[
\begin{matrix}
0&0&0&0&0&0&0&0&0&0 \\
0&y_{45}&-y_4&-y_5&0&0&0&0&0&0 \\
0&-y_3&y_{35}&-y_5&0&0&0&0&0&0 \\
0&-y_3&-y_4&y_{34}&0&0&0&0&0&0 \\
0&0&0&0&y_{45}&-y_4&-y_5&0&0&0 \\
0&0&0&0&-y_3&y_{35}&-y_5&0&0&0 \\
0&0&0&0&-y_3&-y_4&y_{34}&0&0&0 \\
0&0&0&0&0&0&0&y_{345}&0&0 \\
0&0&0&0&0&0&0&0&y_{345}&0 \\
0&0&0&0&0&0&0&0&0&y_{345}
\end{matrix}
\right].
\]

A calculation with the projection $p:\sA^\bul(\cG) \to \sA^\bul(\cS)$
yields the induced endomorphism
$\omega^\bul(\cS',\cS)=\omega^\bul(S,r):\sA^\bul(\cS)\to \sA^\bul(\cS)$, given explicitly 
by $\omega^0(\cS',\cS)=0$, $\omega^1(\cS',\cS)=\tilde\omega^1(S,r)$, 
and
\[
\omega^2(\cS',\cS)= 
\left[
\begin{matrix}
y_{45}&-y_4&-y_5&0&0&0\\
-y_3&y_{35}&-y_5&0&0&0\\
-y_3&-y_4&y_{34}&0&0&0\\
0&0&0&y_{45}&-y_4&-y_5\\
0&0&0&-y_3&y_{35}&-y_5\\
0&0&0&-y_3&-y_4&y_{34}
\end{matrix}
\right].
\]

Weights $\bl=(\la_1,\la_2,\la_3,\la_4,\la_5)$ are nonresonant for type $\cS$ if
\[
\la_1,\  \la_2,\  \la_3,\  \la_4,\  \la_5,\  \la_6,\  \la_{135},\ 
\la_{245},\  \la_{126},\  \la_{346}
\notin \Z_{\ge 0},
\]
where $\la_J=\sum_{j\in J}\la_j$ and $\la_6 = -\la_{[5]}$.
The
$\beta\textbf{nbc}$ basis for $H^2(A^\bul(\cS),a_\bl)=H^2(\sfM;\LL)$ is 
$\{\eta_{2,4},\eta_{2,5}\}$, 
where $\eta_{2,j}=(\la_2 a_2 +\la_4 a_4+\la_5 a_5)\la_j a_j$, 
see \cite{FaT1}.  The  projection map  
$\rho^2:A^2(\cS) \twoheadrightarrow H^2(A^\bul(\cS),a_\bl)$ is given by
\[
\rho^2(a_{i,j})=
\begin{cases}
(\la_{3}\eta_{2,4}-\la_1(\eta_{2,4}+\eta_{2,5}))/(\la_1\la_3\la_{135})
&\text{if $\{i,j\}=\{1,3\}$,}\\
-\eta_{2,4}/(\la_1\la_4)
& \text{if $\{i,j\}=\{1,4\}$,}\\
(\la_{15}\eta_{2,4}+\la_1(\eta_{2,4}+\eta_{2,5}))/(\la_1\la_5\la_{135})
&\text{if $\{i,j\}=\{1,5\}$,}\\
-(\eta_{2,4}+\eta_{2,5})/(\la_2\la_3) &\text{if $\{i,j\}=\{2,3\}$,}\\
(\la_{24}\eta_{2,4}+\la_4\eta_{2,5})/(\la_2\la_4\la_{245})&\text{if $\{i,j\}=\{2,4\}$,}\\
(\la_{5}\eta_{2,4}+\la_{25}\eta_{2,5})/(\la_2\la_5\la_{245})&\text{if $\{i,j\}=\{2,5\}$.}
\end{cases}
\]

A calculation with the endomorphism
$\omega_\bl^2(\cS',\cS)=
\left.\omega^2(\cS',\cS)\right|_{\b{y}\mapsto\bl}$
and this projection yields
\[
\Omega_{\cC}^2(\cS',\cS)=
\left[\begin{matrix}
\la_3+\la_4+\la_5 & 0 \\ 0 & \la_3+\la_4+\la_5
\end{matrix}\right].
\]

A collection of weights $\bl$ is resonant for type $\cS$ if
$\la_1 = \la_4$, $\la_2 = \la_3$, $\la_5 = \la_6$, and 
$\la_1 + \la_2 + \la_5 = 0$.  
Let $\bl$ be a collection of nontrivial, resonant weights.
Then $\la_1\neq 0$ or $\la_2\neq 0$.
For such weights, one can 
check that $a_1-a_2-a_3+a_4 \in A^1(\cS)$ represents a basis for $H^1(A^\bul(\cS),
a_\bl)$, and that $\dim H^2(A^\bul(\cS),a_\bl)=3$.  By Theorem~\ref{thm:evalues}, 
the spectrum of the Gauss-Manin endomorphism $\Omega^q_\cC(\cS',\cS)$ is 
contained in 
$\{0,\la_{345}\}$, provided $\la_{345} \neq 0$.  However, the resonance conditions above 
imply that $\la_{345}=0$.  Accordingly, one can check directly that the endomorphism 
$\Omega^1_{\cC}(\cS',\cS):H^1(A^\bul(\cS),a_\bl)
\to H^1(A^\bul(\cS), a_\bl)$
induced by $\omega^\bul_\bl(\cS',\cS)=
\left.\omega^\bul(\cS',\cS)\right|_{\b{y}\mapsto\bl}$
is trivial.  One can also show that, 
for an appropriate choice of basis for 
$H^2(A^\bul(\cS),a_\bl)$, the projection 
$A^2(\cS) \twoheadrightarrow H^2(A^\bul(\cS),a_\bl)$ 
has matrix
\[
\left[
\begin{matrix}
\la_1 & -\la_1 & \la_2 \\
-\la_2 & \la_2 & \la_2 \\
0 & 0 & \la_2 \\
\la_1 & -\la_1 & -\la_1 \\
\la_1 & \la_2 & \la_2 \\
\la_1 & 0 & 0
\end{matrix}
\right],
\]
and that the endomorphism 
$\Omega^2_{\cC}(\cS',\cS):H^2(A^\bul(\cS),a_\bl)
\to H^2(A^\bul(\cS), a_\bl)$
induced by $\omega^\bul_\bl(\cS',\cS)$
is trivial as well.

\end{document}